\documentclass{article}

\usepackage[utf8]{inputenc}
\usepackage{geometry}[margin=1in]
\usepackage{amsmath, amsthm, amssymb, amsfonts}
\usepackage{xcolor}

\usepackage{stackengine}
\usepackage{bm}
\usepackage[mathscr]{euscript}
\usepackage{enumitem}

\usepackage[backend=bibtex8,autocite = superscript]{biblatex}
\bibliography{literature} 
\usepackage[colorlinks,allcolors=blue]{hyperref}
\let\cite=\supercite

\usepackage{algorithm}
\usepackage[noend]{algpseudocode}

\usepackage{graphicx}
\usepackage[labelfont=bf]{caption}
\usepackage[format=hang]{subcaption}
\usepackage{pstricks, pst-node}

\usepackage{scalerel}

\usepackage{booktabs}
\usepackage{multirow}
\usepackage{float}
\usepackage{colortbl}
\usepackage{makecell}
\makeatletter
\@ifundefined{c@rownum}{%
  \let\c@rownum\rownum
}{}

\newcommand{\RR}{\mathbb{R}}

\DeclareMathOperator*{\argmax}{arg\,max}

\newcommand{\GMSalgoHeader}[1]{\Statex \underline{\textsc{#1}}}

\title{Advancing Geometry with AI:\\Multi-agent Generation of Polytopes}
\date{\today}
\author{Grzegorz \'Swirszcz \and Adam Zsolt Wagner\thanks{Equal contributors. Authors listed alphabetically by last name.} \and  Geordie Williamson\footnotemark[1] \and Sam Blackwell\thanks{Equal contributors. Authors listed alphabetically by last name.} \and Bogdan Georgiev\footnotemark[2] \and Alex Davies\thanks{Equal contributors. Authors listed alphabetically by last name.} \and Ali Eslami\footnotemark[3] \and S\'{e}bastien Racani\`{e}re\footnotemark[3] \and Theophane Weber\footnotemark[3] \ \and Pushmeet Kohli}

\begin{document}

\maketitle

\begin{abstract}
Polytopes are one of the most primitive concepts underlying geometry. They are defined as intersections of a finite number of half-spaces, i.e. the subsets of Euclidean spaces which can be defined by a finite set of linear inequalities. Discovery and study of polytopes with complex structures provides a means of advancing scientific knowledge. Construction of polytopes with specific extremal structure is very difficult and time-consuming. Having an automated tool for the generation of such extremal examples is therefore of great value. We present an Artificial Intelligence system capable of generating novel polytopes with very high complexity, whose abilities we demonstrate in three different and challenging scenarios: the Hirsch Conjecture, the k-neighbourly problem and the longest monotone paths problem. For each of these three problems the system was able to generate novel examples, which match or surpass the best previously known bounds. Our main focus was the Hirsch Conjecture, which had remained an open problem for over 50 years. The highly parallel A.I. system presented in this paper was able to generate millions of examples, with many of them surpassing best known previous results and possessing properties not present in the earlier human-constructed examples. For comparison, it took leading human experts over 50 years to handcraft the first example of a polytope exceeding the bound conjectured by Hirsch, and in the decade since humans were able to construct only a scarce few families of such counterexample polytopes. With the adoption of computer-aided methods, the creation of new examples of mathematical objects stops being a domain reserved only for human expertise. Advances in A.I. provide mathematicians with yet another powerful tool in advancing mathematical knowledge. The results presented demonstrate that A.I. is capable of addressing problems in geometry recognized as extremely hard, and also to produce extremal examples different in nature from the ones constructed by humans.
\end{abstract}

\section{Introduction}
The discovery of interesting examples is a critical component of mathematical progress. The icosahedron, whose discovery is often attributed to the ancient Greeks, is an early example of exceptional structure in mathematics, and features prominently in Euclid’s elements. More recently, the Lorenz attractor, Conway’s game of life and the Mandelbrot set are all examples of 20th century discoveries that gave birth to new fields of mathematical inquiry and new understanding of our world. According to mathematician Gil Kalai ``The methods for coming up with useful examples in mathematics \dots\ are even less clear than the methods for proving mathematical statements.”\cite{ZieglerExtreme}

In this work we focus on \textit{polytopes}, which are higher dimensional generalizations of familiar geometric objects like the triangle, square, cube and tetrahedron. As well as being fundamental objects in geometry, polytopes are also of great importance in optimization. Despite their importance, many mysteries abound concerning their structure. In particular, it can be a very difficult task to construct examples of polytopes with specific extremal structure in high dimensions.

In this paper we propose a new A.I. guided algorithm for the generation of interesting polytopes. A key idea is to regard the generation of polytopes as a one-player game, in which new polytopes are built by starting with a randomly generated list of vertices and then repeatedly selecting and moving a vertex. In high dimensions there are an enormous number of ways to modify a polytope by moving a vertex. This makes any exhaustive search computationally intractable, because the number of possible moves is too large. Moreover, amongst all possible modifications of a given polytope the vast majority of modifications are uninteresting. For this reason, we train an A.I. agent to predict moves in this game which will lead to interesting polytopes.

We demonstrate the power of our approach on three classical problems in the theory of high-dimensional polytopes. The first is the 50-year-old Hirsch conjecture, which asks for a bound on the diameter of a polytope in terms of its dimension and number of facets. In a breakthrough 2012 paper~\cite{santos2012counterexample}, Santos constructed the first counterexample to the Hirsch conjecture, and since his work only two more examples and one family have been discovered~\cite{matschke2015width}. However, each involved many months of human effort and closely related questions still remain unsolved~\cite{polymath, santos2013recent}. Our system is able to generate millions of new counter-examples which appear fundamentally different to the previously known human constructions, and it found a non-Hirsch polytope in only 19 dimensions, improving on the previous bound (20 dimensions). We also test our system on two other problems in polytope theory: the neighbourly polytopes problem and the longest monotone paths problem. In all three cases, we were able to generate new examples that match or surpass the best known bounds.

\section{Background on polytopes}

\begin{figure}[ht]
\centering
\includegraphics[width=14cm]{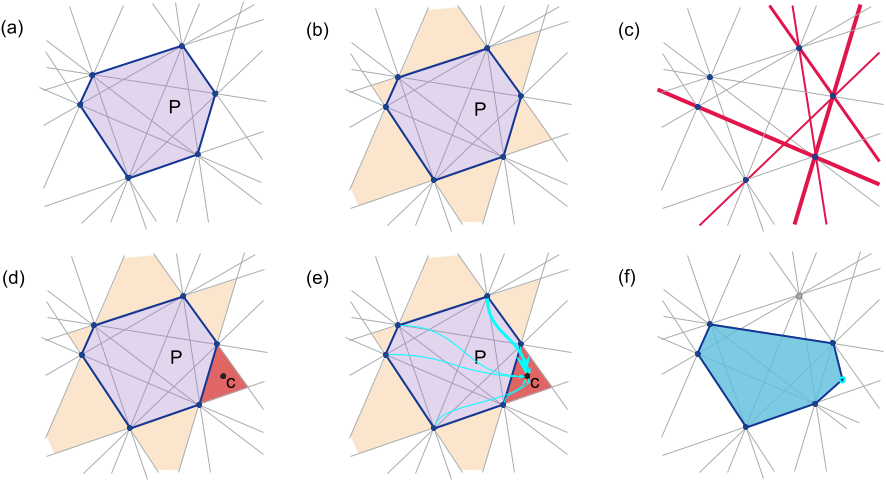}
\caption{An outline of the Hopper algorithm. (a) A polytope $P$ and all the hyperplanes determined by its vertices. (b) Of these regions, adding a point is only possible in the light orange shaded regions. In each of these regions, adding a point yields a combinatorially non-equivalent polytope. (c) After computation of all hyperplanes (red), the neural network determines a probability distribution on these hyperplanes. (d) A region (red) described by a subset of the hyperplanes is selected using the neural network, then a point $c \in R$ is chosen (deterministically). (e) All polytopes obtained by ``hopping" a vertex of $P$ to $c$ are returned, degenerate ones are discarded, and non-degenerate ones are evaluated. (f) One of the resulting polytopes. (This 2$d$ picture obscures some subtleties present in higher dimension. Most importantly, individual modifications can have an enormous effect on the combinatorial type of the polytope. Additionally, the selection of the region $R$ in step (c) is more complicated than this picture makes out. Finally, in higher dimensions the number of regions $R$ yielding potential hops is enormous.)}
\label{fig::magic_table}
\end{figure}

A polytope is the convex hull of finitely many points inside a finite-dimensional vector space $\mathbb{R}^d$ (see Figure \ref{fig::magic_table}(a)). We will always assume that the vertices are not contained in a hyperplane, and that our polytopes are simplicial, meaning that each $(d-1)$-dimensional face is a simplex, i.e.~the convex hull of exactly $d$ vertices. For example, the tetrahedron, octahedron and icosahedron are simplicial, but the cube and regular dodecahedron are not.

Most questions in the combinatorics of polytopes only concern their combinatorial structure and not the explicit position of their vertices. A polytope has faces of dimensions 0, 1, …, $d-1$ and their incidences determine the combinatorial type of a polytope. We will often regard polytopes as the same if they are of the same combinatorial type. Small perturbations of the vertices of a simplicial polytope result in combinatorially equivalent polytopes.\cite{ZieglerPolytopes, de2010triangulations}

In general, the computation of the faces of a polytope from its vertices presents both conceptual and computational challenges as the number of faces of polytopes can be very large. The enumeration of the facets of a polytope given only its vertices is a close relative of the linear programming problem, ubiquitous in many optimization problems.\cite{fukuda1994combinatorial,dantzig1963linear}

Consider a polytope $P$ with $n$ vertices and suppose we wish to modify it by either moving one of its vertices yielding a (combinatorially) different polytope with $n$ vertices, adding a vertex to $P$ to create a new polytope with $n+1$ vertices or deleting a vertex, creating a polytope with $n-1$ vertices.  Several difficulties arise:
\begin{enumerate}
    \item The first two modifications often lead to a polytope with a different number of vertices than desired. This can happen due to some vertices ending up within a convex hull of the others as a result of the modification.
    \item Among the modifications which preserve the number of vertices, many lead to combinatorially equivalent polytopes. The regions leading to combinatorially equivalent polytopes are themselves polytopes, whose faces are given by hyperplanes through the vertices of $P$ (see Figure~\ref{fig::magic_table}(b)). It is computationally expensive to compute all these regions, and there are an enormous number in high dimensional space (e.g. $5\cdot 10^{23}$ for a 30-vertex polytope in 5 dimensions).
    \item For the regions constituting valid moves, moves leading to extremal polytopes are very rare if they exist at all. Thus they are unlikely to be found via random search.
\end{enumerate}
Several innovations are necessary in order to overcome these difficulties, which we now explain.

\section{The Hopper algorithm}
\label{subs::hopper_algorithm}
Throughout we work with polytopes $P$, described via a list of vertices. We fix a property $\rho$ which we would like to achieve (e.g., “facet-ridge graph has diameter $n$”). We also specify several fitness functions $f_i$ which are continuous approximations of how close $P$ is to achieving $\rho$.

\begin{figure}[ht]
\centering
\includegraphics[width=12cm]{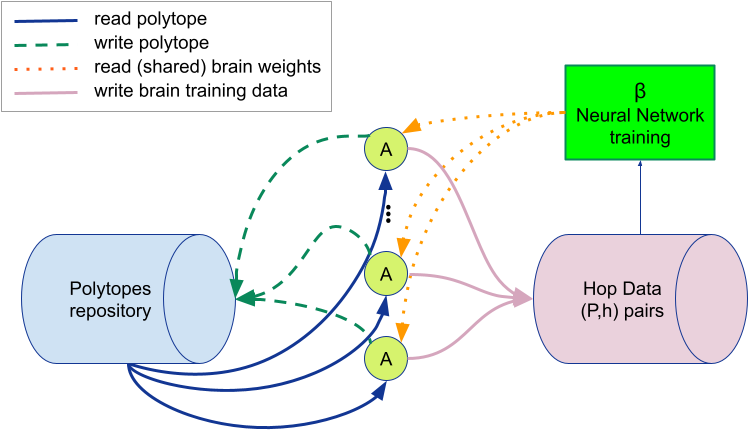}
\caption{\label{fig::multi_agent_architecture} Multi-agent architecture. The agents (A) read and write polytopes from the shared polytope repository. Simultaneously pairs of polytopes and hyperplanes are written to the shared data repository. The pairs are labeled positively if the hyperplane was part of a boundary of a successful hop region and negatively otherwise.}
\end{figure}

While machine learning methods for polytopes have been considered in the past to predict various standard properties~\cite{bao2021polytopes, coates2023machine,  he2023machine, he2023machine1} and databases of polytopes have been built before~\cite{paffenholz2017polydb}, our approach is rather different. The global structure of the algorithm is population refinement. A collection of agents is working together on a shared repository of examples, collectively trying to improve the fitness functions $f_i$, while each of the individual agents has its individual fitness function fixed. The agents operate using a work-stealing strategy. Each agent picks a polytope from a shared repository and attempts to modify it. If an improvement is achieved, the resulting polytope is written back to the repository (see Figure \ref{fig::multi_agent_architecture}). Each inserted polytope has an “importance score” attached, which is a combination of its fitness function values, its history of improvements, etc. If the size of the repository exceeds a fixed maximum size, some stored polytopes are discarded (see Methods, \ref{subs::shared_repository}).

The subtlety in our algorithm arises in the operation of each agent. As we discussed above, the task of moving, adding or deleting a point to create a polytope with a set of desirable properties, is a difficult one. In order to overcome these difficulties we propose the Hopper algorithm, where each agent randomly selects moves and uses a neural network to predict how good these moves are. This is used to select which moves to action. Given a $d$-dimensional polytope $P$, the Hopper algorithm first computes all hyperplanes $\mathcal{H}$ through subsets of size $d$ of vertices of $P$. Then, these hyperplanes are each fed through the neural network, which takes as input the vertices of $P$ as well as the equation of the given hyperplane. The neural network then outputs a probability that this hyperplane is the boundary of a region which constitutes a good move.

Guided by these probabilities, the Hopper algorithm chooses a region $A$ whose boundary is described via hyperplanes belonging to $\mathcal{H}$. (Hyperplanes assigned a high probability by the neural network are more likely to feature in the boundary of $A$ -- see Figure \ref{fig::magic_table}(c) and Methods.) It then computes the center of the largest ball which can be inscribed inside $A$, and by “hopping” each vertex of $P$ to $c$ creates $n$ new polytope candidates\footnote{In the ``addition" mode it also creates a candidate with $n+1$ vertices by adding $c$ as a new vertex.} (see Figure \ref{fig::magic_table}(d,e)). Degenerate examples are discarded, and the fitness function values of the resulting polytopes are evaluated.

\begin{algorithm}[H]\caption{Hopper algorithm - individual agent}
\label{alg::individual_agent}
\begin{center}
\begin{algorithmic}[1]
\Require Polytope $P \in \RR^d$, represented as a $n \times d$ matrix $V$ of vertices. Fitness function $f$ of the agent. Likelihood function $\beta$ \footnotemark.
\Ensure The given vertices span $P$ properly - no vertex is inside a convex hull of other vertices
\GMSalgoHeader{Compute hyperplanes}
\State $\mathcal{H} = \emptyset$ 
\For{ all $\{ v_{1}, \ldots v_{d}\} \subseteq V$ sets of distinct $d$ vertices}
\State $\mathcal{H}.{\rm append}(h)$, $h$ - equation of a hyperplane containing $\{ v_{1}, \ldots v_{d}\}$
\State Compute probability distribution $\pi$ on hyperplanes in $\mathcal{H}$ by computing likelihoods of each $h$ using likelihood function $\beta$ (see~\ref{subs::hopper_brain}) and normalizing them by dividing by their sum. 
\EndFor
\GMSalgoHeader{Construct hop target $c$}
\State $\mathcal{A} = \emptyset$ 
\State Sample $d+1$ equations $h_1, \ldots h_{d+1}$ from $\mathcal{H}$ according to $\pi$ and change them to inequalities $\chi_1, \ldots, \chi_{d+1}$, choosing directions of the inequalities such that the feasible set is nonempty
If the feasible set is not bounded or the shape violates geometric constrains, discard $\chi_1, \ldots, \chi_{d+1}$ and sample again
\State $\mathcal{A} \leftarrow \chi_1, \ldots, \chi_{d+1}$
\State Find $(c, r)$ - the center and radius of the largest ball inscribed in the polytope defined by $\mathcal{A}$
\While{exists $h \in \mathcal{H}$ such that ${\rm dist}(h, c) < 0.8 * r$}
\State $\mathcal{A}.{\rm append}(h)$
\State Recompute $(c, r)$
\EndWhile
\GMSalgoHeader{Construct list of candidates}
\State Q = []
\For{$v \in V$}
\State Construct polytope $P'$ by replacing $v$ with $c$
\If{$P'$ satisfies admissibility criteria}
\State $Q.{\rm append}(P')$
\EndIf
\EndFor
\GMSalgoHeader{Evaluate candidates}
\For{$\varrho \in Q$}
\If{$f(\varrho) < f(P)$}
\State $Q.{\rm delete}(\varrho)$
\EndIf
\EndFor
\State \Return $\argmax\limits_f Q$
\end{algorithmic}
\end{center}
\end{algorithm}
\footnotetext{Likelihood function $\beta$ is an output of a neural network being trained in an online fashion. For more details see Methods~\ref{supp::hopper_brain_architecture}.}

\subsection{Using a neural network for selecting actions}
\label{subs::hopper_brain}
Initially the actions are chosen uniformly at random. As the work progresses, more and more knowledge can be exploited. The agents are equipped with a (shared) ``brain" $\beta$ that allows this to happen. The architecture of that brain -- a transformer-based neural network -- is described in detail in~\ref{supp::hopper_brain_architecture}. The brain takes as an input a polytope and a hyperplane pair $(P, h)$ and outputs a floating point number representing a likelihood of hyperplane $h$ being part of a boundary of a successful hop region (see Algorithm~\ref{alg::individual_agent}). The training happens in an online mode -- the (labelled) data from each agent is written to a shared data repository (``Hop Data" silo on Figure~\ref{fig::multi_agent_architecture}) from where it is sampled by the training process.
In reinforcement learning terminology, the brain $\beta$ corresponds to the policy, which we are improving over time by learning from our experiences.


\section{Polytope discovery results}
\label{sec::results}

\subsection{Hirsch conjecture -- finding non-Hirsch polytopes}
\label{subs::hirsch_main}
A fundamental problem in linear programming is to bound the number of iterations the simplex method needs to take to solve an optimization problem that has $d$ variables and $n$ inequalities~\cite{borgwardt2012simplex}.  To find a lower bound for this quantity, we would like to know what the largest possible distance between any two vertices can be, in a $d$-dimensional polytope with $n$ facets. The celebrated Hirsch conjecture asserts that the diameter of a $d$-dimensional convex polytope with $n$ facets is at most $n-d$. To date, the best upper bound known is quasipolynomial~\cite{sukegawa2019asymptotically, todd2014improved, kalai1992quasi}. 

A \emph{prismatoid} is a special class of polytopes introduced by Santos. A prismatoid is a polytope with two parallel facets (the \emph{base facets}), which together contain all vertices. The \emph{width} of a prismatoid is the distance, in the facet-ridge graph, of its base facets. Santos explained how to start with a $d$-dimensional prismatoid with
$n$ vertices of width $\ge w$ and produce another prismatoid with $2(n-d)$ vertices in
dimension $n-d$ of diameter $\ge w+n-2d$. Thus, the Hirsch conjecture
implies that prismatoids in dimension $d$ have width at most $d$. 

This concept of prismatoids turned out to be very fruitful in attacking the Hirsch conjecture. Indeed, every counterexample humans have found for the Hirsch conjecture originated from a 5-dimensional prismatoid. We followed this idea, and use the Hopper algorithm to generate over a million new 5-dimensional prismatoids with width 6. The smallest counterexample we found has only 24 vertices, beating the state of the art of 25 vertices~\cite{matschke2015width}. This prismatoid has 12 points in both base facets and width 6, see section~\ref{supp::data}. By the discussion above it corresponds to a non-Hirsch polytope of dimension 19 (24-5), the smallest known counterexample so far.

\subsubsection{New features}
The constructions found by the Hopper algorithm seem to be of a genuinely different nature to known human constructions. It has been observed by Matschke, Santos and Weibel that all known counterexamples to the Hirsch conjecture may be explained by a simple observation: one may construct a bipartite directed graph from any prismatoid via a simple geometric procedure, and if this digraph contains no directed cycles of length 2, then the polytope is non-Hirsch (\cite{matschke2015width}, Proposition 2.3). We have checked that this explains all human constructions of non-Hirsch polytopes. We checked a sample of 100 of the non-Hirsch polytopes constructed by the Hopper algorithm, and none of them can be explained in this way. This indicates that the polytopes constructed by the Hopper algorithm are fundamentally non-human, potentially providing new lines of investigation in the field.

\begin{figure}[ht] 
\centering
\begin{subfigure}[t]{0.34\textwidth}
\includegraphics[width=0.9\textwidth]{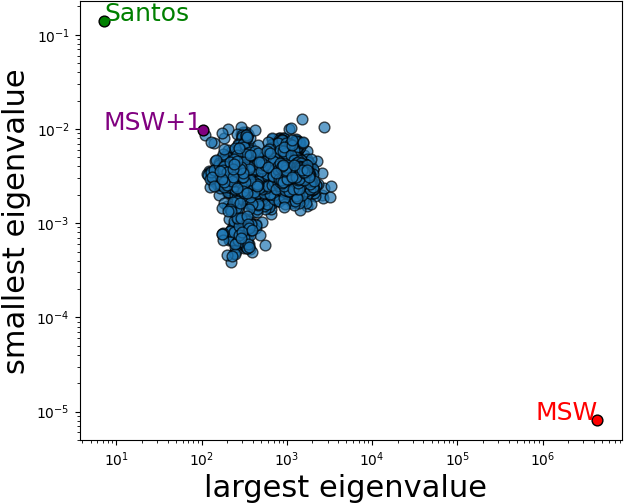}
\caption{PCA scales of our constructions versus human constructions.}
\label{fig::pca}
\end{subfigure}\hfill
\begin{subfigure}[t]{0.30\textwidth}
\includegraphics[width=0.9\textwidth]{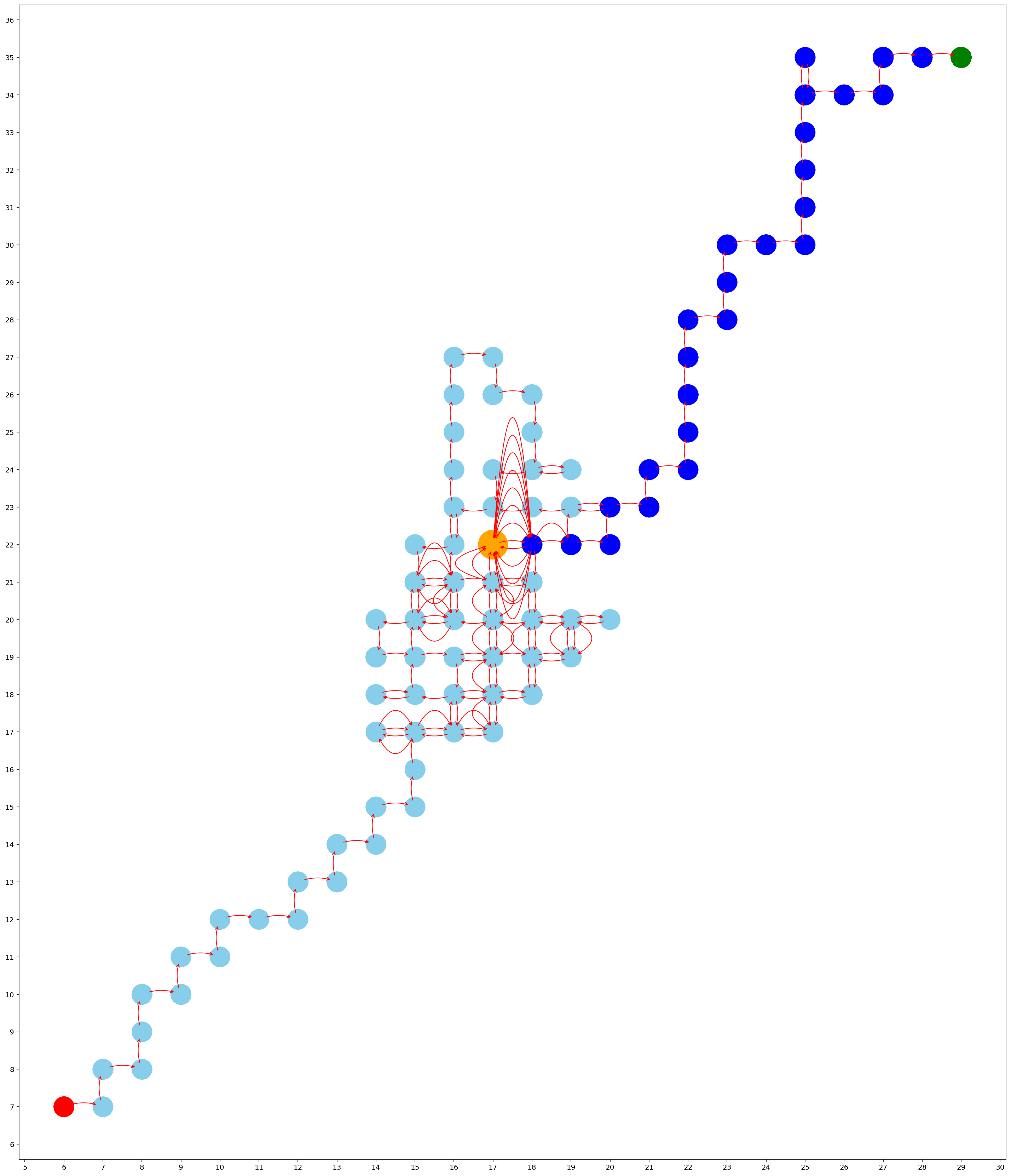}
\caption{Prismatoid generation: 937, defect: 95,  final number of vertices: $29\times35$, ascension\protect\footnotemark at $17\times22$ vertices.}
\end{subfigure}\hfill
\begin{subfigure}[t]{0.345\textwidth}
\includegraphics[width=0.9\textwidth]{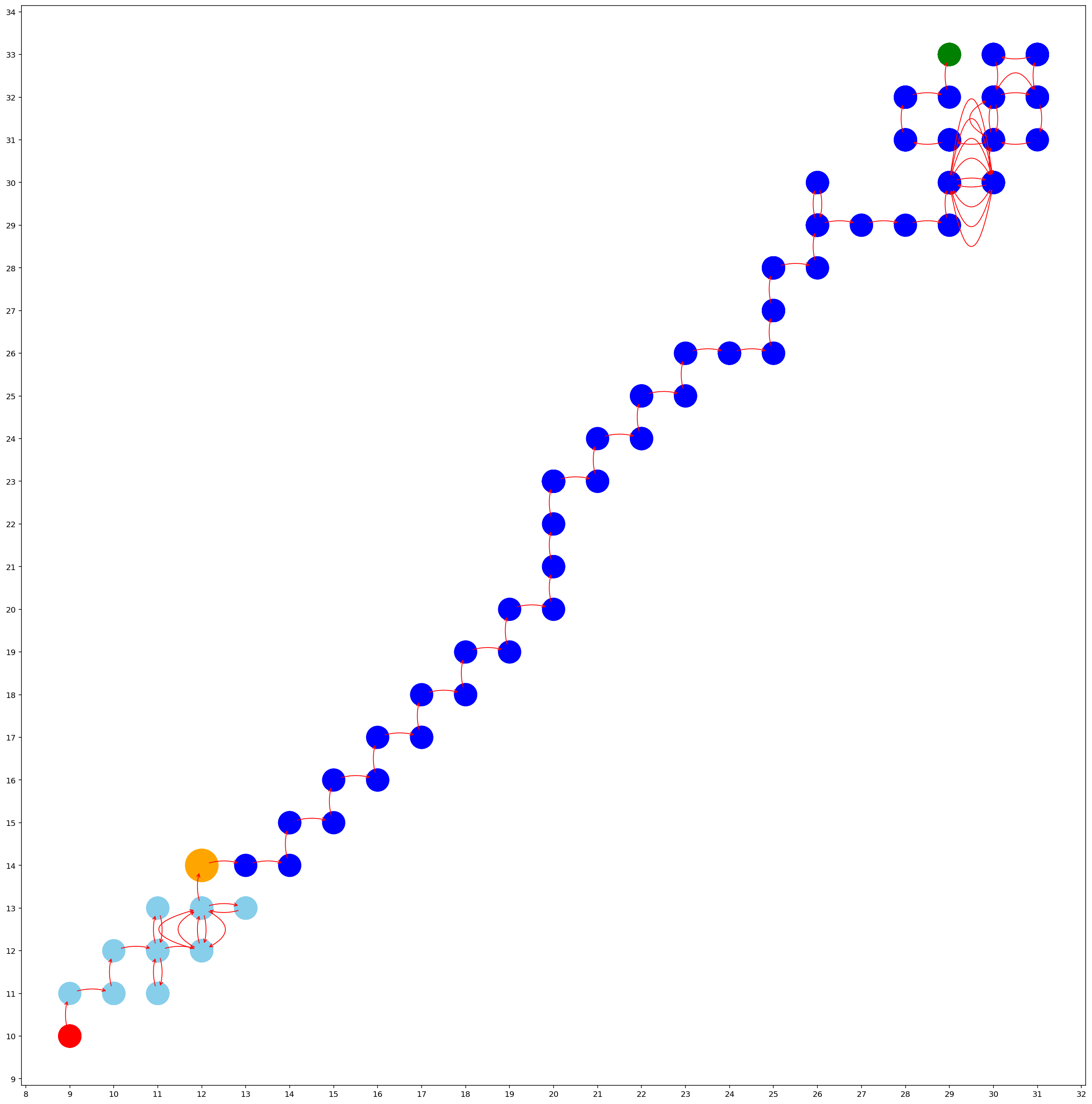}
\addtocounter{footnote}{-1}\addtocounter{Hfootnote}{-1} 
\caption{Prismatoid generation: 850, defect: 11, final number of vertices: 29x33, ascension\protect\footnotemark at $12\times14$ vertices.}
\end{subfigure}
\addtocounter{footnote}{-1}\addtocounter{Hfootnote}{-1} 
\caption{(a) The different scales of the constructions given by the Hopper algorithm, versus human examples. We normalize the bottom base facet and record the PCA eigenvalues of the top base facet, then plot the largest eigenvalue versus the smallest eigenvalue at log scale. The green, red, and purple dots correspond to Santos's original construction, the 25 vertex constructions in~\cite{matschke2015width}, and the first member of the family of constructions in~\cite{matschke2015width} (all other members of this family have approximately the same position in this plot). The blue dots are the constructions found by the Hopper algorithm. We can see the examples found by the Hopper algorithm feature several orders of magnitude of difference between the smallest and largest PCA eigenvalues. 
(b) Evolution of a prismatoid. The axes represent numbers of vertices. The cyan dots represent the states before the prismatoid increases its width to 6, the blue ones correspond to the states after. The orange dot marks the ``ascension" \protect\footnotemark~to the non-Hirsch realm of width 6 prismatoids. 
(c) Evolution of another prismatoid. This prismatoid has a low {\it defect } (see~\ref{supp::heuristics}) equal to $11$.} 
\end{figure}

We also see new features when considering the global geometry of prismatoids yielding counterexamples to the Hirsch conjecture. Examples constructed by Maschke, Santos and Weibel relied on polytopes with drastically different coordinate scales -- their coordinates range from $10^{-3}$ to $10^{9}$.
One way of analysing this is to apply a linear transformation to the bottom base facet, so that all of its PCA eigenvalues are 1, and then analyse the PCA eigenvalues of the top base facet. The largest and smallest PCA eigenvalues of the prismatoids discovered by Santos and Maschke-Santos-Weibel are presented in Figure~\ref{fig::pca}. It is striking that all examples found by the Hopper algorithm involve multiple orders of magnitude of difference between the lowest and highest PCA eigenvalues.  It is also interesting to note that the infinite family found by Maschke, Santos and Weibel has coordinates in a range very close to those discovered by the Hopper algorithm.

\footnotetext{The Hopper algorithm starts off with a population of polytopes satisfying the Hirsch condition ${\rm diameter} \le n-d $. By means of sequence of hops some of the polytopes become non-Hirsch polytopes. We called the hop that attains that the ``ascension". After becoming non-Hirsch the polytope continues to be worked on, hence the blue part of the trajectory.}

\subsection{The longest monotone path problem}
\label{subs::monotone_paths}
Given a polytope, how long is the longest `monotone path' in it? That is, the longest sequence of vertices, each consecutive pair connected by an edge in the polytope, so that a certain linear functional is always increasing. 

Among polytopes with $n$ facets, in dimension $d$, denote the longest possible such monotone path by $f(n, d)$.   Many general constructions are given in~\cite{klee1972good,  amenta1999deformed, pfeifle2004monotone}. We can use the length of the longest path along which the $x$ coordinate always increases as a fitness function. To be able to generate polytopes with exactly $n$ facets we worked in the dual space. Using the Hopper algorithm, we obtained the following lower bounds, beating the state-of-the-art for each pair of parameters that we tried.

\begin{center}
\begin{tabular}{ c|c|c } 
 $(n,d)$ & previously known bounds & new bounds from the Hopper algorithm \\ \hline
 $(9,5)$ & $1 \leq f(9, 5) \leq 30$ & $f(9,5) = \textbf{30}$ \\ 
 $(10,5)$ & $40 \leq f(10, 5) \leq 42$ & $\textbf{41} \leq f(10,5) \leq 42$  \\ 
 $(11,5)$ & $32 \leq f(11, 5) \leq 56$ & $\textbf{55} \leq f(11,5) \leq 56$\\
 $(9,6)$ & $27 \leq f(9, 6) \leq 29$ & $f(9,6) = \textbf{29}$
\end{tabular}
\end{center}

\subsection{The neighbourly polytopes problem}
\label{subs::k_neighbourly_polytopes}
 A polytope is called \emph{neighbourly} if every set of at most $ \lfloor d/2 \rfloor$ vertices forms a face. For any $(n,d)$ pair, there exists a unique $d$-dimensional polytope with $n$ vertices called the \emph{cyclic polytope}, and every cyclic polytope is neighbourly.
Motzkin~\cite{gale1963neighborly} conjectured that all neighborly polytopes are  equivalent to cyclic polytopes. This turned out to be false. It is now known that there exist polytopes in $d\geq 4$ dimensions with $d+4$ or more vertices that are neighbourly but not cyclic. There has been a lot of work on finding new classes of neighbourly polytopes  via ingenius human constructions, see e.g.~\cite{shemer1982neighborly, shemer1984techniques, devyatov2011neighbourly}. To calculate the number of faces of various dimensions, we used the algorithm presented in~\cite{fukuda1994combinatorial}. The fitness function was then the largest $k$ for which every set of at most $k$ vertices forms a face, plus the proportion of $k+1$-sets of vertices that form a face. Using this fitness function, the Hopper algorithm was easily able to find many previously unknown neighbourly polytopes that are not cyclic, for several pairs $(n,d)$ ranging from $(10,6)$ to $(13,8)$. 

\section{Conclusion}
\label{sec::conclusion}
We presented a new A.I. guided algorithm that when trained from scratch, discovered new examples of polytopes that are novel and exceed the best known results on three classical problems in geometry. 
An important strength of the Hopper algorithm is its flexibility -- it is a general purpose algorithm that can easily be adapted to find polytopes with a wide variety of extremal properties. It supports the optimization of complex non-differentiable rewards, and when there is not a single good choice for a reward function, it is capable of working with a larger set of approximate reward functions. This flexibility enables the user to apply the Hopper algorithm to complex problems where they may not have a full understanding of the precise reward function they want to optimize for. Another advantage of our approach is that the constructions it produces are different from constructions found by humans, which can lead to further insights and discoveries. 

The automatic discovery of polytopes with never before seen properties has the potential to advance mathematics. What took humans many years, can be discovered by the Hopper algorithm much faster. Our approach demonstrates the power of A.I.and its ability to assist in making new mathematical discoveries in a field of mathematics which has resisted effective computer search due to its complexity.
 
\newpage

\appendix

\section{Methods}

\subsection{Preliminaries}
Let $P$ be a convex polytope. The connection graph of $P$ has nodes and edges
corresponding to the vertices and edges of $P$. The diameter of $P$ is the maximum distance between any two points in the graph of $P$.  For some problems, it is  useful to dualize. To $P$ we may also associate its facet-ridge graph which
has nodes given by the facets (=codimension 1 faces) of $P$ and edges corresponding
to ridges (=codimension 2 faces) of $P$. The facet-ridge graph of $P$ agrees with the
(vertex-edge) connection graph of the dual polytope. 
To work with polytopes computationally, while there are other specialized tools available~\cite{gawrilow2000polymake}, we relied on the scipy and pycddlib packages (see~\ref{supp::stability}).

\subsection{Our approach}
Numerous successes have been obtained by applying Reinforcement Learning (RL) to various problems from playing games to making scientific breakthroughs, see for example \cite{Fawzi2022, davies2021advancing,  wagner2021constructions}. The concept of RL is very broad and general, hence our approach could be presented within this framework. Nevertheless, in the traditional Reinforcement Learning the main focus is on the policy improvement. In our case, the main pillars of the solution are the multi-agent knowledge distillation and the concept of multi-objective optimization described in \ref{supp::multi-objective}.

\subsection{Memoization}
\label{subs::memoization}
Memoization means storing the results of a computed expressions so that they can be subsequently retrieved without repeating the computation. This technique is used frequently in the code. One prominent example would be in recomputing the planes in Algorithm~\ref{alg::individual_agent}. After the hop is executed not all hyperplanes change, and only the ones that changed are being recomputed in order to improve performance.
\subsection{Batching}
When computing the probability distributions over hyperplanes (see Algorithm~\ref{alg::individual_agent}) the neural network is applied to batches of inputs rather than one individual hyperplane at a time, for performance reasons. This is a very standard practice in the field of Neural Networks.
\subsection{"Safety fuses"}
In most {\it while} loops break clauses were introduced based on timeouts and/or numbers of repetitions. This prevented accidental freezing or unnecessary slowdowns of the agents. One prominent use of said technique was the construction of a hop region - the break was executed if the addition of subsequent hyperplanes was taking too long, if the number of attempts to construct the initial simplex was too high etc.
\subsection{Numerical stability of computations}
\label{supp::stability}
The objects we are working with are usually very complicated polytopes. Hence it is crucial to not get false results as an outcome of numerical roundoff errors. The CDD library offers the ability to work with exact fractions, at a cost of performance/memory utilization. The Hopper works in the dual mode. Many pre-selection procedures rely on the representation using float numbers. Final evaluation of candidates always uses the precise mode.
\subsection{Multi-objective optimization}
\label{supp::multi-objective}
\begin{figure}[ht]
\centering
\begin{subfigure}[t]{0.48\textwidth}
\includegraphics[width=\textwidth]{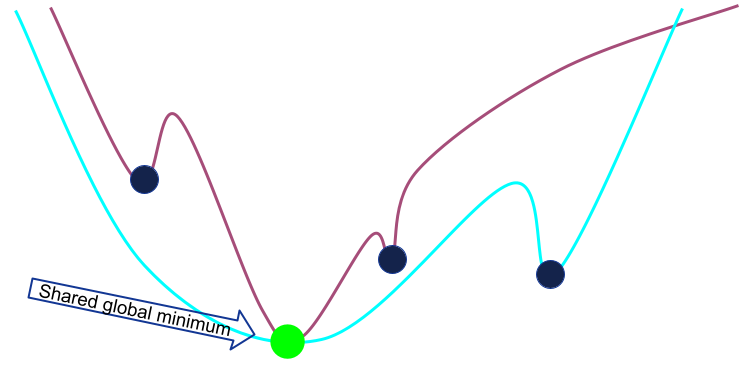}
\caption{Success mode 
\label{fig::multi_objective_success}}
\end{subfigure}\hfill
\begin{subfigure}[t]{0.48\textwidth}
\includegraphics[width=\textwidth]{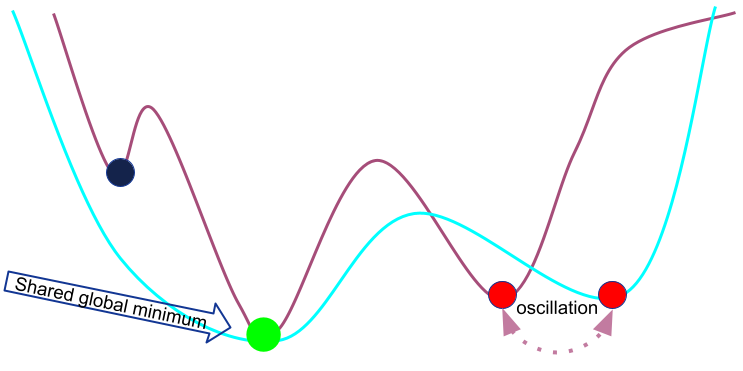}
\caption{Failure - oscillating between two suboptimal local minima
\label{fig::multi_objective_failure}}
\end{subfigure}
\caption{\label{fig::multi_objective_training} Multi-objective training. Two fitness functions share the same global minimum corresponding to a desired state. Alternating between the objectives can lead to bypassing a barrier of local minima.}
\end{figure}

If for an optimization problem it is possible to have more than one fitness function, it can be often used to the advantage of the optimization process. One of the most notorious problems of optimization is being stuck in "bad regions", such as local minima or plateaus. If the two fitness functions share a global minimum corresponding to the desired solution - but differ from one another, one can hope for their "bad regions" to not overlap, as in Figure~\ref{fig::multi_objective_success}. Then a savvy switching procedure between the fitness function can make us successfully bypass the local minima and arrive at the global minimum. The procedure has its inherent failure modes, as illustrated in Figure~\ref{fig::multi_objective_failure}.
In finding counterexamples to Hirsch conjecture a collection of 10 different fitness functions was successfully used. Those fitness functions were a {\it defect} (see~\ref{supp::heuristics}) and a collection of formulae counting the shortest and longest paths connecting the top and bottom base facet.

\subsection{Heuristics}
\label{supp::heuristics}
In the Hirsch conjecture, the primary metric we want to optimize is the width of the prismatoid, i.e.~the distance between the two base facets in the facet-ridge graph. To make this fitness metric more continuous, a natural quantity to define is the average distance between neighbors of the bottom base facet and neighbors of the top base facet. This ``average width'' was a useful general concept, as it enabled us to have a more continuous measure of progress. To use multi-objective optimization (see~\ref{supp::multi-objective}) we worked with several variants of this quantity. 
In case of the input being a prismatoid, a very important\footnote{In our opinion - also elegant.} one among those is a {\it defect}, defined as follows:
Let $l_{{\rm min}}$ be the length of the shortest path connecting top to the bottom. The {\it defect} of a prismatoid is the total number of paths of length $l_{{\rm min}}$ from a vertex in the top face to a vertex in the bottom face.

 \subsection{Adding and subtracting vertices}
 \label{subs::adding_and_deleting}

For applications to the Hirsch conjecture, it was necessary to extend the agent to allow the addition and subtraction of vertices.\footnote{The first new counterexample to Hirsch conjecture found by Hopper algorithm was operating in a rigid $12\times13$ regime, with no addition or deletion allowed. Enabling the addition and deletion changed the expected time to finding the first example from weeks to hours, and enabled to produce millions of them.}  However, we found that doing this naively led to degenerate behaviour: either the agents removed vertices until we were left with simplices, or the agent kept adding vertices until computations became infeasible. This degenerate behaviour was avoided by the addition of preventive rules.

 The Hopper can operate in two modes. In the {\it rigid} mode the Hopper is allowed only the modifications by taking a vertex and moving it somewhere else. In this mode the total number of vertices is preserved. In a {\it flexible} mode the Hopper is allowed two more modifications: deleting a vertex and adding a vertex. The latter modification needs an explanation. Recall that the basic modification (see~\ref{subs::hopper_algorithm}) consists of choosing a vertex $v_i$ and a new vertex $\bar{v}$ inside a region $V_{\iota}$ and then replacing the vertex $v_i$ with the vertex $\bar{v}$. The operation of addition means simply adding the vertex $\bar{v}$ without removing any existing vertices, hence increasing the number of vertices. In case such modification would not lead to the increase of vertices due to one or more of existing vertices becoming the interior points or any other disallowed geometric condition the modification is discarded. Deletion of vertices is done randomly. For addition, we use the neural network to select a region in a way similar to that describe above, see~\ref{subs::hopper_brain}.
 We discovered that simply enabling the flexible mode leads to trivial behaviours. In almost every case the vertices are either removed until the polytope becomes a simplex, or are being added with the number of vertices tending towards infinity. To remedy the tendency to converge to one of those trivial regimes a heuristic was invented. It was called an uptick and downtick rule. Uptick means that the addition of a vertex is only allowed, if the total number of the shortest paths (compare~\ref{supp::heuristics}) decreases as a result of it. This means that the addition did something "interesting" to the geometry, as the natural tendency of the addition of a vertex is to increase the number of vertices, and hence - the total number of connections. The downtick rule acts in a dual manner, using a slightly technical concept of "longest paths".\\
 In the "production" runs of the Hirsch conjecture experiments (\ref{subs::hirsch_main}) used the {\it flexible} mode, while the monotone paths experiments (\ref{subs::monotone_paths}) and k-neighbourly polytopes experiments (\ref{subs::k_neighbourly_polytopes}) used the {\it rigid} mode.

\subsection{Avoiding geometric degeneration}
\label{supp::geometry}
As the agent modifies the polytope, it might propose changes leading to undesired geometric properties of the polytope. Examples of such properties would be very large coordinates of vertices (escape to infinity), very small angles between the facets, extreme disproportions between sizes in different directions etc. Some of those issues can be remedied by simply applying an affine transformation. A heuristic maintaining a "canonical scaling" based on Principal Components Analysis is invoked after every update. Other issues are handled by a number of heuristics which do not admit hops that lead to flatness, disproportions etc.

\subsection{``Hopper brain" architecture}
\label{supp::hopper_brain_architecture}
The hopper brain takes as an input a pair - a $d$-dimensional prismatoid and a hyperplane determined by $d$ vertices of the prismatoid. The prismatoid and the hyperplane data is first preprocessed independently, then is is stacked together and fed into the final Transformer block, followed by a Feed Forward block. The approach follows the pattern from the~\cite{vaswani2017attention} and uses multi-headed attention and positional encoding. The details are as follows.
\subsubsection{Preprocessing the prismatoid}
\begin{itemize}
    \item Each base facet is separately processed by a 2-headed attention block with 3 layers and a key size $d+3$
    \item Each base facet is further fed through GELU~\cite{hendrycks2016gelu} activated Feed Forward block with two linear layers
    \item Top and bottom base facets are embelished with a $[1, 0]$ and $[0, 1]$ positional encodings respectively and stacked into a full prismatoid matrix
    \item Full prismatoid matrix is fed into a a 2-headed attention block with 4 layers and a key size $3d+1$
    \item The rows corresponding to the top base facet are further positionally encoded with $[0,0, 1,0,1,0,0]$
    \item The rows corresponding to the bottom base facet are further positionally encoded with $[0,0,0,1,1,0,0]$
\end{itemize}

\subsubsection{Preprocessing the hyperplane}
A (signed) hyperplane-constraint is encoded as a single vector, corresponding to the inequality $Ax \le b$.
A positional encoding token ${\rm DECK}$ is created, equal to $[1,0]$ if the hyperplane is a part of the hop in the top base facet and $[0,1]$ if it is a part of the hop in the bottom base facet.
\begin{itemize}
    \item Hyperplane is processed by a gelu activated Feed Forward block with 3 layers
    \item Positionally encode processed hyperplane with $[{\rm DECK}, 0,0,0,1,0]$
    \item Positionally encode raw hyperplane with $[{\rm DECK}, 0,0,0,0,1]$
\end{itemize}
\subsubsection{Joint stage}
\begin{itemize}
    \item Stack vertically the processed hyperplane, the raw hyperplane and the processed prismatoid together with their positional encoding\footnote{Thanks to the proper choice of sizes of layers, all dimensions match.}
    \item Apply 2-headed attention with 4 layers and key size $4 d  + 1$
    \item Process the output of the main attention block by $4$-layer gelu activated Feed Forward block, decreasing the dimensionality to $3$
\end{itemize}

\subsubsection{Main loss}
The main loss is computed as a cross-entropy (see for example~\cite{zhang2018generalized}) loss between the output of the final Feed Forward block and a $1$-hot $3$-dimensional vector indicating the ground truth of one of the scenarios: $[1,0,0]$ - the hyperplane was a part of a partition cell that led to a successful hop, $[0,1,0]$ if the hyperplane was a part of a partition cell that got rejected due to geometrical infeasibility (see~\ref{supp::heuristics}) and $[0,0,1]$ if the hyperplane was a part of a partition cell which was geometrically feasible, but did not lead to a successful hop. A hop is considered successful, if it lead to an improvement in a heuristic and a subsequent increase in the heuristic value.

\subsubsection{Auxiliary loss}
The pre-processed prismatoid is fed into a $4$-layer gelu activated Feed Forward block, decreasing the dimensionality to $128$. The auxiliary loss is then computed as a cross entropy between the output and a $128$-dimensional 1-hot vector denoting the ground truth value of the {\it defect} (see~\ref{supp::heuristics}) of the prismatoid, with the vector with 1 at the last position indicates any defect greater or equal to $127$.
The purpose of the auxiliary loss is to "teach the brain some geometry". This follows the philosophy of AI and the paradigm of a transfer of knowledge. The prediction of the defect n itself is not important to us, but it improves the quality of the model.
\subsubsection{Total loss and the training}
The total loss is the sum of the main loss and the auxiliary loss. The model was implemented using the JAX library and trained using Adam~\cite{adam} algorithm in Optax optimizer. 

\begin{figure}[ht]
\centering
\subfloat[a][]{\includegraphics[width=\textwidth]{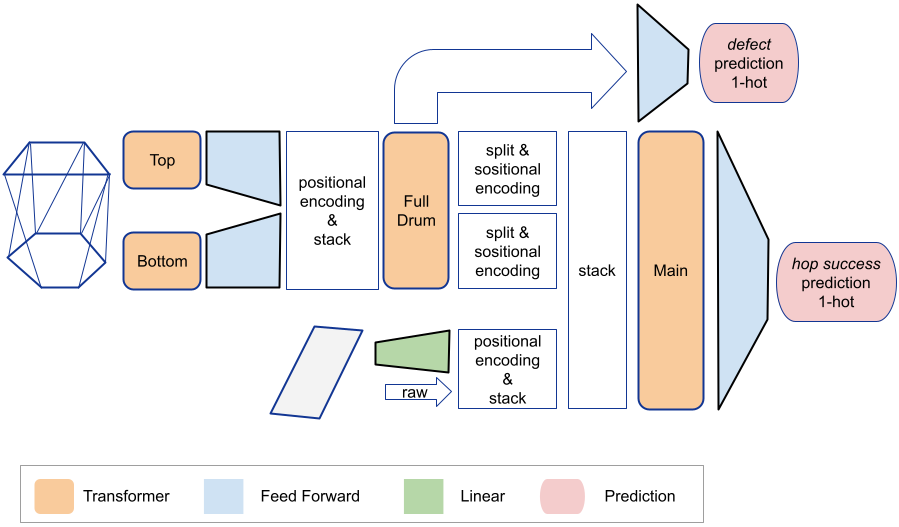}}
\\
\subfloat[b][]
{\includegraphics[width=\textwidth, height=\textwidth/3]{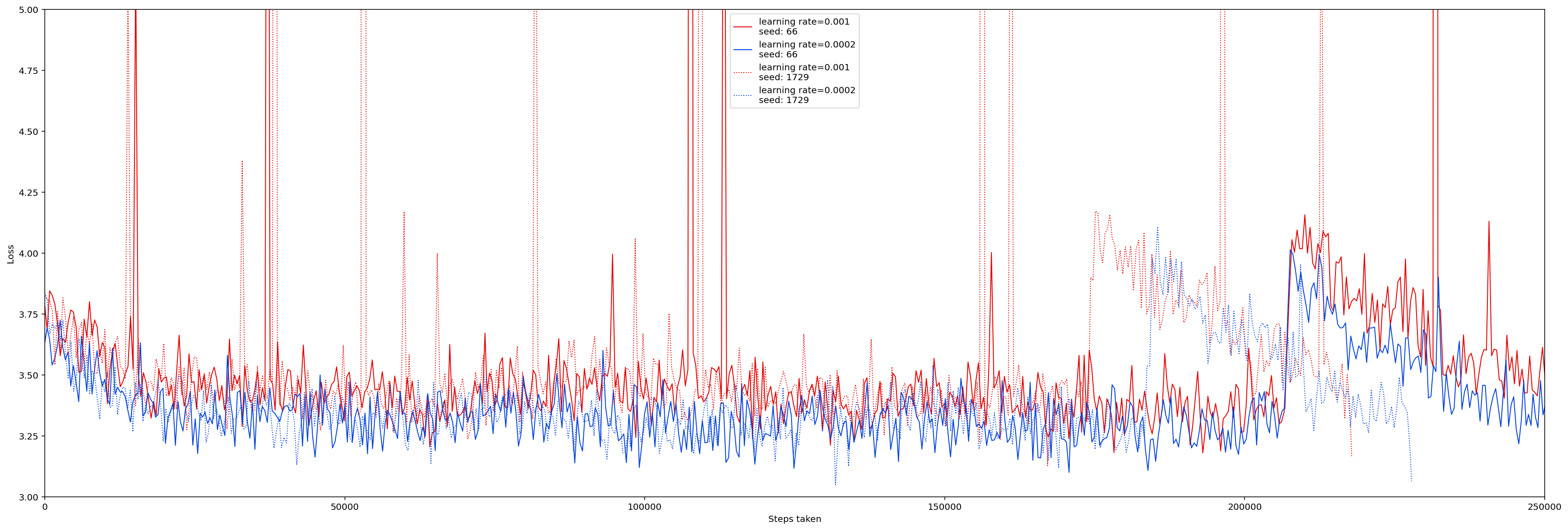}}
\caption{\label{fig::hopper_architecture} {\bf (a)} The architecture of the Hopper brain. The input is a pair -- a polytope and a hyperplane. The outputs are: the main prediction, of whether the input hyperplane was in a boundary of some successful hop region, and the auxiliary prediction -- the {\it defect} (see~\ref{supp::heuristics}) of the input prismatoid. The "arrow block" highlights the fact that the hyperplane data leg does not participate in the prediction of the defect. {\bf (b)} Training loss of the Hopper brain. As opposed to more traditional scenarios the training is happening in an online fashion. As the population of polytopes in the repository becomes more refined the model can be expected to lose its validity at time. Such change of the regime can be well seen around step $2.77M$.}
\end{figure}


\subsection{Shared repository \& Reverb}
\label{subs::shared_repository}
In Subsection~\ref{subs::hopper_algorithm} the multiagent architecture of Hopper swarm is explained. It uses two data repositories, one for the polytope examples, and one for the purpose of training the shared "consciousness" of the agents. In order for such repositories to be of practical use a lot of performance demands had to be satisfied. The speed, reliability, well implemented addition/removal protocols etc. A Google DeepMind solution called Reverb~\cite{cassirer2021reverb} was chosen for the task. Reverb is a replay buffer designed to support large scale Reinforcement Learning experiments and has proven to be a great fit for the needs of the Hopper.

\label{supp::shared_repository}

\subsection{Data}
\label{supp::data}

The coordinates of the 24-vertex prismatoid of width 6 are as follows.
\vspace{1cm}

{
\footnotesize
\begin{tabular}{rrrrr}
\hline
 -0.002&  1.302& -0.072& -0.277&  1\\
        0.028& -1.854& -0.031& -1.509&  1\\
       -0.031& -2.024&  0.062& -0.648&  1\\
        0.083& -0.544& -0.048&  1.771&  1\\
       -0.024&  1.084&  0.082&  0.168&  1\\
        0.047&  0.108&  0.046& -2.533&  1\\
       -0.028&  1.317&  0.002& -1.173&  1\\
        0.013& -1.755&  0.051&  0.348&  1\\
       -0.037& -2.020& -0.066& -0.423&  1\\
        0.001& -1.320&  0.018& -2.118&  1\\
       -0.086& -1.040&  0.019&  1.657&  1\\
       -0.061&  1.190& -0.029& -1.518&  1\\
        1.751&  0.035& -0.574& -0.029& -1\\
       -0.821& -0.017&  2.210&  0.216& -1\\
       -1.696&  0.033& -0.938& -0.046& -1\\
        1.531& -0.024& -0.932& -0.032& -1\\
       -1.743& -0.041&  0.440& -0.008& -1\\
        1.344&  0.034&  1.333&  0.125& -1\\
        1.696& -0.045&  0.587& -0.025& -1\\
        0.670&  0.005&  1.961& -0.077& -1\\
        0.692&  0.020& -1.894& -0.134& -1\\
       -0.156& -0.002& -1.900&  0.540& -1\\
       -1.573& -0.029& -1.286& -0.171& -1\\
       -1.693&  0.052&  0.995& -0.011& -1\\
\hline
\end{tabular}
}

This prismatoid has a {\it defect} (see~\ref{supp::heuristics}) equal to $64$, which is the smallest value we observed for a $12 \times 12$ example of a non-Hirsch prismatoid.

\printbibliography

@misc{vaswani2017attention,
      title={Attention Is All You Need}, 
      author={Ashish Vaswani and Noam Shazeer and Niki Parmar and Jakob Uszkoreit and Llion Jones and Aidan N. Gomez and Lukasz Kaiser and Illia Polosukhin},
      year={2017},
      eprint={1706.03762},
      archivePrefix={arXiv},
      primaryClass={cs.CL}
}

@misc{cassirer2021reverb,
      title={Reverb: A Framework For Experience Replay},
      author={Albin Cassirer and Gabriel Barth-Maron and Eugene Brevdo and Sabela Ramos and Toby Boyd and Thibault Sottiaux and Manuel Kroiss},
      year={2021},
      eprint={2102.04736},
      archivePrefix={arXiv},
      primaryClass={cs.LG}
}

@article{adam,
  title={Adam: A Method for Stochastic Optimization},
  author={Diederik P. Kingma and Jimmy Ba},
  journal={CoRR},
  year={2015},
  volume={abs/1412.6980}
}

@Article{Fawzi2022,
author={Fawzi, Alhussein
and Balog, Matej
and Huang, Aja
and Hubert, Thomas
and Romera-Paredes, Bernardino
and Barekatain, Mohammadamin
and Novikov, Alexander
and R. Ruiz, Francisco J.
and Schrittwieser, Julian
and Swirszcz, Grzegorz
and Silver, David
and Hassabis, Demis
and Kohli, Pushmeet},
title={Discovering faster matrix multiplication algorithms with reinforcement learning},
journal={Nature},
year={2022},
month={Oct},
day={01},
volume={610},
number={7930},
pages={47-53},
issn={1476-4687},
doi={10.1038/s41586-022-05172-4},
url={https://doi.org/10.1038/s41586-022-05172-4}
}

@article{hendrycks2016gelu,
  title={Learning to reason in large theories without imitation},
  author={Hendrycks, Dan and Gimpel, Kevin},
  journal={arXiv:1606.08415},
  year={2016}
}

@article{santos2012counterexample,
  title={A counterexample to the {H}irsch conjecture},
  author={Santos, Francisco},
  journal={Annals of mathematics},
  volume={176},
  issue={1},
  pages={383--412},
  year={2012}
}

@article{wagner2021constructions,
  title={Constructions in combinatorics via neural networks},
  author={Wagner, Adam Zsolt},
  journal={arXiv preprint arXiv:2104.14516},
  year={2021}
}

@inproceedings{gale1963neighborly,
  title={Neighborly and cyclic polytopes},
  author={Gale, David},
  booktitle={Proc. Sympos. Pure Math},
  volume={7},
  pages={225--232},
  year={1963}
}

@incollection{shemer1984techniques,
  title={Techniques for investigating neighborly polytopes},
  author={Shemer, Ido},
  booktitle={North-Holland Mathematics Studies},
  volume={87},
  pages={283--292},
  year={1984},
  publisher={Elsevier}
}

@article{devyatov2011neighbourly,
  title={Neighbourly polytopes with few vertices},
  author={Devyatov, Rostislav A},
  journal={Sbornik: Mathematics},
  volume={202},
  number={10},
  pages={1441},
  year={2011},
  publisher={IOP Publishing}
}

@article{shemer1982neighborly,
  title={Neighborly polytopes},
  author={Shemer, Ido},
  journal={Israel Journal of Mathematics},
  volume={43},
  pages={291--314},
  year={1982},
  publisher={Springer}
}

@article{klee1972good,
  title={How good is the simplex algorithm},
  author={Klee, Victor and Minty, George J},
  journal={Inequalities},
  volume={3},
  number={3},
  pages={159--175},
  year={1972},
  publisher={New York}
}

@article{pfeifle2004monotone,
  title={On the monotone upper bound problem},
  author={Pfeifle, Julian and Ziegler, G{\"u}nter M},
  journal={Experimental Mathematics},
  volume={13},
  number={1},
  pages={1--11},
  year={2004},
  publisher={Taylor \& Francis}
}

@article{amenta1999deformed,
  title={Deformed products and maximal shadows of polytopes},
  author={Amenta, Nina and Ziegler, Gunter M},
  journal={Contemporary Mathematics},
  volume={223},
  pages={57--90},
  year={1999},
  publisher={Providence, RI: American Mathematical Society}
}

@article{fukuda1994combinatorial,
  title={Combinatorial face enumeration in convex polytopes},
  author={Fukuda, Komei and Rosta, Vera},
  journal={Computational Geometry},
  volume={4},
  number={4},
  pages={191--198},
  year={1994},
  publisher={Elsevier}
}

@article{matschke2015width,
  title={The width of five-dimensional prismatoids},
  author={Matschke, Benjamin and Santos, Francisco and Weibel, Christophe},
  journal={Proceedings of the London Mathematical Society},
  volume={110},
  number={3},
  pages={647--672},
  year={2015},
  publisher={Oxford University Press}
}

@article{polymath,
title={Polymath 3: Polynomial Hirsch Conjecture},
author={Gil Kalai (coordinator)},
url={http://gilkalai.wordpress.com/2010/09/29/polymath-3-polynomial-hirsch-conjecture},
year={2010}
}

@article{santos2013recent,
  title={Recent progress on the combinatorial diameter of polytopes and simplicial complexes},
  author={Santos, Francisco},
  journal={Top},
  volume={21},
  number={3},
  pages={426--460},
  year={2013},
  publisher={Springer}
}

@incollection {ZieglerExtreme,
    AUTHOR = {Ziegler, G\"{u}nter M.},
     TITLE = {Convex polytopes: extremal constructions and {$f$}-vector
              shapes},
 BOOKTITLE = {Geometric combinatorics},
    SERIES = {IAS/Park City Math. Ser.},
    VOLUME = {13},
     PAGES = {617--691},
 PUBLISHER = {Amer. Math. Soc., Providence, RI},
      YEAR = {2007},
      ISBN = {978-0-8218-3736-8; 0-8218-3736-2},
   MRCLASS = {52B05 (52A10)},
  MRNUMBER = {2383133},
       DOI = {10.1090/pcms/013/10},
       URL = {https://doi.org/10.1090/pcms/013/10},
}

@book{dantzig1963linear,
  title={Linear programming and extensions},
  author={Dantzig, George},
  year={1963},
  publisher={Princeton university press}
}

@book {ZieglerPolytopes,
    AUTHOR = {Ziegler, G\"{u}nter M.},
     TITLE = {Lectures on polytopes},
    SERIES = {Graduate Texts in Mathematics},
    VOLUME = {152},
 PUBLISHER = {Springer-Verlag, New York},
      YEAR = {1995},
     PAGES = {x+370},
      ISBN = {0-387-94365-X},
   MRCLASS = {52Bxx},
  MRNUMBER = {1311028},
MRREVIEWER = {Margaret\ M.\ Bayer},
       DOI = {10.1007/978-1-4613-8431-1},
       URL = {https://doi.org/10.1007/978-1-4613-8431-1},
}

@article{sukegawa2019asymptotically,
  title={An asymptotically improved upper bound on the diameter of polyhedra},
  author={Sukegawa, Noriyoshi},
  journal={Discrete \& Computational Geometry},
  volume={62},
  pages={690--699},
  year={2019},
  publisher={Springer}
}

@article{todd2014improved,
  title={An improved Kalai--Kleitman bound for the diameter of a polyhedron},
  author={Todd, Michael J},
  journal={SIAM Journal on Discrete Mathematics},
  volume={28},
  number={4},
  pages={1944--1947},
  year={2014},
  publisher={SIAM}
}

@article{kalai1992quasi,
  title={A quasi-polynomial bound for the diameter of graphs of polyhedra},
  author={Kalai, Gil and Kleitman, Daniel J},
  journal={Bulletin of the American Mathematical Society},
  volume={26},
  number={2},
  pages={315--316},
  year={1992}
}

@article{davies2021advancing,
  title={Advancing mathematics by guiding human intuition with AI},
  author={Davies, Alex and Veli{\v{c}}kovi{\'c}, Petar and Buesing, Lars and Blackwell, Sam and Zheng, Daniel and Toma{\v{s}}ev, Nenad and Tanburn, Richard and Battaglia, Peter and Blundell, Charles and Juh{\'a}sz, Andr{\'a}s and others},
  journal={Nature},
  volume={600},
  number={7887},
  pages={70--74},
  year={2021},
  publisher={Nature Publishing Group UK London}
}

@incollection{coates2023machine,
  title={Machine learning: The dimension of a polytope},
  author={Coates, Tom and Hofscheier, Johannes and Kasprzyk, Alexander M},
  booktitle={Machine Learning in Pure Mathematics and Thoeretical Physics},
  pages={85--104},
  year={2023},
  publisher={World Scientific}
}

@article{bao2021polytopes,
  title={Polytopes and machine learning},
  author={Bao, Jiakang and He, Yang-Hui and Hirst, Edward and Hofscheier, Johannes and Kasprzyk, Alexander and Majumder, Suvajit},
  journal={arXiv preprint arXiv:2109.09602},
  year={2021}
}

@article{paffenholz2017polydb,
  title={polyDB: a database for polytopes and related objects},
  author={Paffenholz, Andreas},
  journal={Algorithmic and experimental methods in algebra, geometry, and number theory},
  pages={533--547},
  year={2017},
  publisher={Springer}
}

@article{he2023machine,
  title={Machine-learning mathematical structures},
  author={He, Yang-Hui},
  journal={International Journal of Data Science in the Mathematical Sciences},
  volume={1},
  number={01},
  pages={23--47},
  year={2023},
  publisher={World Scientific}
}

@book{de2010triangulations,
  title={Triangulations: structures for algorithms and applications},
  author={De Loera, Jes{\'u}s and Rambau, J{\"o}rg and Santos, Francisco},
  volume={25},
  year={2010},
  publisher={Springer Science \& Business Media}
}

@inproceedings{gawrilow2000polymake,
  title={Polymake: a framework for analyzing convex polytopes},
  author={Gawrilow, Ewgenij and Joswig, Michael},
  booktitle={Polytopes—combinatorics and computation},
  pages={43--73},
  year={2000},
  organization={Springer}
}

@article{he2023machine1,
  title={Machine Learning in Physics and Geometry},
  author={He, Yang-Hui and Heyes, Elli and Hirst, Edward},
  journal={arXiv preprint arXiv:2303.12626},
  year={2023}
}

@book{borgwardt2012simplex,
  title={The simplex method: a probabilistic analysis},
  author={Borgwardt, Karl Heinz},
  volume={1},
  year={2012},
  publisher={Springer Science \& Business Media}
}

@article{zhang2018generalized,
  title={Generalized cross entropy loss for training deep neural networks with noisy labels},
  author={Zhang, Zhilu and Sabuncu, Mert},
  journal={Advances in neural information processing systems},
  volume={31},
  year={2018}
}

\end{document}